\newcommand{\bull}{\vrule height .9ex width .8ex depth -.1ex}
 \newcommand{\ppp}{\hfill $\bull$ }
 \newcommand{\hhh}{{\rm Hess}\,}
 \documentclass[11pt]{article}
 \author{ M.-L. Labbi}
   \title{On $(2k)$-Minimal Submanifolds}
   \date{}
   \usepackage{amsmath}
   \usepackage{amscd}
\newtheorem{theorem}{Theorem}[section]
\newtheorem{corollary}[theorem]{Corollary}
\newtheorem{lemma}[theorem]{Lemma}

\newtheorem{proposition}[theorem]{Proposition}
\newtheorem{definition}[theorem]{Definition}
\begin{document}
   \maketitle
   \begin{abstract} Recall that a submanifold of a Riemannian manifold is said to
be minimal if its mean curvature is zero. It is classical that minimal submanifolds 
are the critical points of the volume function.\\
In this paper, we examine  the critical points of the total 
$(2k)$-th Gauss-Bonnet curvature function, called $(2k)$-minimal submanifolds.
We prove that they are characterized by the vanishing of a higher mean curvature,
namely the $(2k+1)$-Gauss-Bonnet curvature.\\
Furthermore, we show that  several properties  of usual minimal submanifolds can be naturally
generalized to
 $(2k)$-minimal submanifolds.

   \end{abstract}
   \par\bigskip\noindent
 {\bf  Mathematics Subject Classification (2000):} 53C40, 53C42.
   \par\medskip\noindent
   {\bf Keywords.}  Generalized minimal submanifolds, generalized Laplacian, Gauss-Bonnet
curvatures.
  
\section{Introduction}
Recall that a submanifold $M$ of a Riemannian manifold $(\tilde M,\tilde g)$ is said to
be minimal if its mean curvature vanishes everywhere. It is classical that minimal submanifolds 
are the critical points of the volume function.\\
In this paper, we consider the critical points of the total 
$(2k)$-th Gauss-Bonnet curvature function, called $(2k)$-minimal submanifolds.\\
These generalize ordinary minimal submanifolds obtained for $k=0$ and  Reilly's $r$-minimal 
hypersurfaces of the Euclidean space when $r=2k$, \cite{Reilly}.\\
We prove that they are characterized by the vanishing of a higher mean curvature,
namely the $(2k+1)$-Gauss-Bonnet curvature. This result generalizes a similar result of 
Reilly obtained for submanifolds of the Euclidean space \cite{Reilly2}.\\
The paper is divided into two parts. In the first part, we first recall useful facts
about some operations on double forms, namely  the exterior product, generalized Hodge star operator,
contraction map and the inner product.\\ 
These tools are used in this paper to provide first
an alternative elegant approach  to symmetric functions and
Newton transformations. And secondly, to provide  a natural introduction to the $k-th$ 
Gauss-Bonnet
curvatures (for $k$ even or odd) and the Einstein-Lovelock tensors.\\
In the second part of this paper, we first prove the first variation formula for the total
$(2k)$-th Gauss-Bonnet curvature function in order to characterize the critical points. Next, 
we prove several facts about $(2k)$-minimal submanifolds which generalize similar properties
about ordinary minimal submanifolds. In particular, we prove that
complex submanifolds of a Kahlerian manifold are always $(2k)$-minimal for all $k$. Also,
we show that compact, irreducible isotropy homogeneous spaces always admit $(2k)$-minimal
immersions in a sphere.\\
The natural extension of Laplace operator that naturally appears in our context is the 
operator
$\ell_{2k}$. Roughly speaking, it is obtained by "contracting" the Hessian by the
$(2k)$-th Einstein-Lovelock
tensor, (recall that the usual Laplacian is just the contraction of the Hessian 
by the metric under consideration).\\
We prove that for a compact manifold these generalized Laplacians are self adjoint
and with zero integral, in fact they can be written as a divergence. Furthermore,
if the metric on the manifold has positive (resp. negative) definite $(2k)$-Einstein-Lovelock
tensor then the operator $\ell_{2k}$ is elliptic and positive (resp. negative) definite. 
In particular,
we obtain a maximum principle for these operators.\\
Finally, we study some properties of $(2k)$-minimal immersions in Euclidean space and spheres.
For example, we prove that an isometric immersion $F: M\rightarrow \Re^{n+p}$ is $(2k)$-minimal
if and only if  the coordinates functions $F_i$ of $F$ are $\ell_{2k}$-harmonic functions.
In particular, we prove that 
there are no non trivial compact $(2k)$-minimal submanifolds in the Euclidean space
with positive definite (or negative definite) $(2k)$-th Einstein-Lovelock tensor.\\

\section{Elementary symmetric functions vs. Gauss-Bonnet curvatures}
\subsection{Double Forms: Algebraic properties}

Let $(V,g)$ be an Euclidean real vector space  of dimension $n$. In the
 following
we shall identify whenever convenient (via their Euclidean structures),
the vector spaces
 with their duals. Let
  $\Lambda V=\bigoplus_{p\geq 0}\Lambda^{p}V$ denotes the exterior algebra
 of  $p$-vectors on $V$. \\
A double form on $V$ of degree $(p,q)$ is  defined to be  a 
 bilinear form $\Lambda^pV\times\Lambda^qV\rightarrow {\bf R}$. Alternatively, it  is
     a multilinear form defined on $V$ which is skew symmetric in the first $p$-arguments and also
     in the last $q$-arguments. If $p=q$ and the bilinear form is symmetric we say that we
 have a symmetric double form.\\
The usual exterior product of $p$-vectors extends in a natural way to double forms of any degree
\cite{Labbidoubleforms}. In particular, the exterior product of two ordinary bilinear forms
coincides with the Kulkarni-Nomizu product. Furthermore, $k$-times the exterior product of a
 symmetric bilinear form  $B$ with itself is
a symmetric double form of order $(k,k)$ and is given by
\begin{equation*}
B^k(x_1 \wedge...\wedge x_k,y_1\wedge...\wedge y_k)=k!\det[B(x_i,y_j)].
\end{equation*}
In particular, for $B=g$, $\frac{g^k}{k!}$ is the canonical inner product on $\Lambda^{k}V$.The former inner
 product
induces a natural inner product of double forms and shall be denoted by $\langle, \rangle$.
The contraction map $c$ on double forms is the adjoint of the exterior 
 multiplication map by the
metric  $g$.\\
Suppose we have chosen  an orientation on the vector space  $V$.
The classical
Hodge star operator  $*:\Lambda^{p}V\rightarrow \Lambda^{n-p}V$ can be extended naturally to
 operate 
on double forms by declaring for a $(p,q)$-double form the following:
\begin{equation*}
 *\omega(.,.)=(-1)^{(p+q)(n-p-q)}\omega(*.,*.).
\end{equation*}
Note that $*\omega$ does not depend on the chosen orientation as the usual Hodge star
 operator is applied twice.
 The so-obtained  operator provides a simple relation between the contraction map $c$ of double forms 
and the multiplication map by the metric:
\begin{equation}
\label{A}
g\omega=*c*\omega\, \, {\rm and}\, c\omega =*g*\omega.
\end{equation}
Furthermore, we have the following properties for all $\omega,\theta\in D^{p,q}$:

\begin{equation}\label{B}
<\omega,\theta>=*(\omega.*\theta)=(-1)^{(p+q)(n-p-q)}*(*\omega.\theta),
\end{equation}
\begin{equation}\label{B+}
**\omega =(-1)^{(p+q)(n-p-q)}\omega.
 \end{equation}
Finally, if $\omega$ is a symmetric $(p,p)$-double form satisfying the first Bianchi
identity then we have \cite{Labbidoubleforms}
\begin{equation}\label{C}
*(\frac{g^{n-p}\omega}{(n-p)!})=\frac{1}{p!}c^p\omega\, \, {\rm and}\, \, 
*(\frac{g^{n-p-1}\omega}{(n-p-1)!})=\frac{c^p\omega}{p!} g -\frac{c^{p-1}\omega}{(p-1)!}.
\end{equation}

\subsection{Elementary Symmetric Functions}
Let $(V,g)$ be an Euclidean space of dimension $n$, $B$ a given  symmetric bilinear form on $V$. 
We denote by $\lambda_1\leq \lambda_2\leq ... \leq \lambda_n$
 the eigenvalues of the
operator corresponding to $B$ via $g$. Let
 $s_k=s_k(\lambda_1,...,\lambda_n)$
be the elementary
symmetric functions for $k=0,1,...,n$, where
 $s_0=1, s_1=\sum_{i=1}^n\lambda_i,...,
 s_n=\lambda_1 ...\lambda_n$.\\
\smallskip\noindent 
 The previous operations on double forms provide an alternative nice 
way to write these invariants as follows.\\
\begin{proposition}
Let $(V,g)$ be an Euclidean space, $B$ a given  symmetric bilinear form.
If $s_k$ denotes the $k$-th elementary symmetric function in the eigenvalues of
 the operator
corresponding to $B$ and $c,*,B^k$ denote respectively the contraction
 map, the generalized Hodge star operator
and the exterior product of $B$ with itself $k$-times then
\begin{equation*}
s_k=\frac{1}{(k!)^2}c^kB^k= \frac{1}{k!(n-k)!}*(g^{n-k}B^k).
\end{equation*}
In particular, the trace and determinant of the operator associated to $B$ via $g$ are given by
\begin{equation*}
s_1={\rm tr}_g  B=*\left\{\frac{g^{n-1}}{(n-1)!}B\right\} \, \, 
{\rm and}\,\,  s_n={\det}_g B=*\frac{B^n}{n!}.
\end{equation*}
\end{proposition}
{\sl Proof}. It is not difficult to see that the eigenvalues of $\frac{B^k}{k!}$ are all  
possible products
$\lambda_{i_1}\lambda_{i_2}...\lambda_{i_k}$ with $i_1<i_2<...<i_k$. From this it is clear that
its complete contraction determines $s_k$ as in the proposition. The second statement 
is a direct application of formula (\ref{C}) above. \ppp \\

\begin{corollary}
Let $A$,$B$ be symmetric bilinear forms,
denote by  $s_i(A),s_j(B)$ the elementary symmetric functions in the eigenvalues of
 the operator
corresponding to $A$ and $B$ respectively via the scalar product $g$. If $A=B+\lambda g$ for some $\lambda\in\Re$, then
for each $k$, $0\leq k\leq n$, we have
\begin{equation*}
s_k(A)=\sum_{i=0}^k\frac{k!(n-i)!}{i!(k-i)!(n-k)!}s_i(B)\lambda^{k-i}.
\end{equation*}
\end{corollary}
{\sl Proof}. Straightforward, just use the binomial theorem and the previous proposition.\ppp \\

\subsection{Newton Transformations}

Associated with the elementary symmetric functions are the so-called Newton transformation 
\cite{Reilly}. We reformulate below their definition in terms of the operations introduced 
above of double forms:
\begin{definition}
For $0\leq k\leq n$, the  $k$-th Newton transformation of a bilinear form $B$ on $(V,g)$
 is defined to be
\begin{equation*}
t_k(B)=* \left\{ \frac{g^{n-k-1}}{(n-k-1)!}
\frac{B^k}{k!}\right\}.
\end{equation*}
For $k=n$, we set $t_n(B)=0$.
\end{definition}
The following   properties of $t_k$ are known \cite{Reilly}:
\begin{proposition}
For each $k$, $0\leq k\leq n$, we have for $B$ and $t_k(B)$ as above
\begin{enumerate}
\item $\langle t_k(B), B\rangle =(k+1)s_{k+1}(B)$. 
This property is equivalent to the celebrated Newton's formula.

\item $t_k(B)=s_k(B)g-\frac{c^{k-1}B^k}{(k-1)!}$.

\item  $ct_k(B)=(n-k)s_k(B).$
\end{enumerate}
\end{proposition}
It is only for the seek to  illustrate the elegance of this new approach to
 symmetric functions and  Newton transformations
that we are proving the previous  proposition below.\\
{\sl Proof.}  To prove the first part, we need just to use  formulas (\ref{B}) and
(\ref{B+}) as follows:
\begin{equation*}
\langle t_k(B), B\rangle =*(\{*t_k(B)\}B)=*\left\{ \frac{g^{n-k-1}}{(n-k-1)!}
\frac{B^{k+1}}{k!}\right\}=(k+1)s_{k+1}(B).
\end{equation*}
The second part  results directly from  formula (\ref{C}).\\
Finally, the third part results from  formulas (\ref{A}) and (\ref{B+}) as follows:
\begin{equation*}
ct_k(B)=*g*t_k(B)=*\left\{ \frac{g^{n-k}}{(n-k-1)!}
\frac{B^k}{k!}\right\}= (n-k)s_k(B).
\end{equation*}
\ppp

\subsection{Gauss-Bonnet curvatures and Einstein-Lovelock tensors}
Let $( M, g)$ be a hypersurface of the $(n+1)$-dimensional Euclidean space.
 The Gauss equation relates the second
fundamental form $B$ of $M$ to its Riemann curvature tensor. Precisely, it states that
$R=1/2B^2$, where of course the product in $B^2=BB$ is the exterior product of double forms. In particular,
$B^{2k}=2^kR^k$ and therefore
the even order symmetric functions in the eigenvalues of $B$ and the corresponding
Newton transformations are  intrinsic invariants of the geometry of the hypersurface
 and are respectively given   by
\begin{equation*}
s_{2k}=\frac{2^k}{[(2k)!]^2}c^{2k}R^{k}=* \left\{ 2^k\frac{g^{n-2k}}{(n-2k)!}
\frac{R^k}{(2k)!}\right\}.
\end{equation*}
\begin{equation*}
t_{2k}=* \left\{ 2^k\frac{g^{n-2k-1}}{(n-2k-1)!}
\frac{R^k}{(2k-1)!}\right\}.
\end{equation*}
The  even order symmetric functions in the eigenvalues of $B$ and the corresponding
Newton transformations are no longer intrinsic for hypersurfaces of arbitrary Riemannian
manifolds. Instead of that, we consider the following  natural intrinsic generalization
 of these curvatures:
\begin{definition} Let $( M, g)$ be an $n$-dimensional Riemannian manifold and 
 let $k$ be a positive integer such that $0\leq 2k\leq n$.
\begin{enumerate}
\item The $(2k)$-th Gauss-Bonnet curvature,   denoted $h_{2k}$,
is the function defined on $M$ by
\begin{equation}\label{D}
h_{2k}={1\over (n-2k)!}*\bigl( g^{n-2k}R^k\bigr).
\end{equation}
\item The $(2k)$-th Einstein-Lovelock tensor, denoted $T_{2k}$, is
 defined by
\begin{equation}\label{D+}
T_{2k}=*{1\over (n-2k-1)!}g^{n-2k-1}R^k.
\end{equation}
\end{enumerate}
If $2k=n$, we set $T_n=o$. For $k=0$ we have $h_0=1$ and $T_0=g$.
\end{definition}
Using formula (\ref{C}) above, these invariants can alternatively be written as
\begin{equation}\label{D++}
h_{2k}=\frac{c^{2k}R^k}{(2k)!}\,\,{\rm and}\,\, T_{2k}=h_{2k}g-\frac{c^{2k-1}R^k}{(2k-1)!}.
\end{equation}
Note that $h_2$ is the half of the usual scalar curvature and $T_2$ is the usual Einstein 
tensor. Recall that if $n$ is even then $h_n$ is up to a constant the Gauss-Bonnet
integrand of $(M,g)$.\\
An important property of these invariants is that the Einstein-Lovelock tensor 
$T_{2k}$ is the gradient of the total $(2k)$-th Gauss-Bonnet curvature seen as a functional
on the space of Riemannian metrics on $M$, see \cite{Labbivariation}.\\

\medskip\noindent
In the special case of a hypersurface of a space form with constant $c$, the invariants
$h_{2k}$ are related to the symmetric functions $s_{2i}$ of the eigenvalues of the second
 fundamental form (which are intrinsic in this special case), as follows:
\begin{proposition}
In a space form of curvature $c$ we have
\begin{equation*}
h_{2k}=\frac{1}{2^k(n-2k)!}\sum_{i=0}^k\frac{k!(n+2i-2k)!}{i!(k-i)!}s_{2k-2i}c^i,
\end{equation*}
and
\begin{equation*}
s_{2k}=\frac{k!}{(n-2k)!}\sum_{i=0}^k(-1)^{k-i}
\frac{2^i(n-2i)!}{i!(k-i)!}h_{2i}c^{k-i}.
\end{equation*}
\end{proposition}
{\sl Proof.} For a hypersurface of a space form with constant $c$, the Gauss equation asserts that
the Riemann curvature tensor of the hypersurface is determined from the second fundamental form $B$
by $R=c\frac{g^2}{2}+\frac{B^2}{2}$. Inserting this in the formulas defining $h_{2k}$ and $s_{2k}$
 we get the desired results.\ppp \\
Similar formulas hold for $T_{2k}$ and $t_{2k}$.

\section{Generalized minimal submanifolds}

Let $(\tilde M,\tilde g)$ be an $(n+p)$-dimensional Riemannian manifold,
 and let $M$ be an
$n$-dimensional  submanifold of $\tilde M$. We shall denote by $g$ the induced metric on $M$.
The purpose of this section is
to characterize those submanifolds (endowed with the induced metric)
 that are critical points of the total Gauss-Bonnet curvature function.\\
\subsection{Gauss-Bonnet curvatures of odd order}
Recall that the Gauss-Bonnet curvatures  $h_{2k}$ of $(M,g)$ are
 intrinsic invariants and are
 defined
by (\ref{D}). We extend the definition of these curvatures to cover odd orders as follows:
\begin{definition}
For a normal vector $N$ at a point $m\in M$ and for $n\geq 2k+1$, we define the
$(2k+1)$ Gauss-Bonnet curvature of the submanifold $(M,g)$ by
\begin{equation}
h_{2k+1}(N)=*(\frac{g^{n-2k-1}}{(n-2k-1)!}R^kB_N).
\end{equation}
For $n=2k$, set $h_{2k+1}(N)=0$.
Where $B$ denotes the vector valued second fundamental form of $M$ and 
$B_N(u,v)=\tilde g(B(u,v),N).$
\end{definition}
The so obtained invariants $h_{2k+1}$ are  normal differential forms
 on $M$ of degree $1$ (duals of  normal vector fields). They are tensorial in $N$.\\
For $k=0$, using (\ref{C}) we get
$$h_1(N)=*(\frac{g^{n-1}}{(n-1)!}B_N)=cB_N.$$
That is $h_1$ is nothing but the usual mean curvature of $M$. Furthermore, for a hypersurface of
the Euclidean space the invariant $h_{2k+1}$ can be seen as a scalar function on $M$ and
$$h_{2k+1}=*\bigl( \frac{g^{n-2k-1}}{(n-2k-1)!}
(\frac{1}{2}B^2)^k B\bigr)=\frac{(2k+1)!}{2^k}s_{2k+1}.$$
That is, up to a constant, the usual $(2k+1)$-mean curvature of the hypersurface $M$.\\
\medskip\noindent 
Using formulas (\ref{B}) and (\ref{B+}), it is straightforward that
\begin{equation}\label{E}
h_{2k+1}(N)=*((*T_{2k})B_N)=\langle T_{2k}, B_N\rangle.
\end{equation}

\subsection{Double Forms: Differential Properties}\label{Diff}
For the seek of completeness, we recall in this paragraph some useful 
differential properties of double forms, for more details see \cite{Kulk, Labbivariation}.\\
Let  $(M,g)$ be a Riemannian manifold of dimension  $n$ and $T_mM$ its tangent space 
at $m\in M$. We denote by   $D^{p,q}$ the vector bundle over $M$ 
whose fiber at $m$ is  the space of all $(p,q)$-double forms on $T_mM$ as in the first part.\\
Note that the previous algebraic properties are still true for the sections of the bundle 
$D^{p,q}$.
\par\medskip\noindent
The second Bianchi map, denoted  $D$, maps  $ D^{p,q}$ into  $ D^{p+1,q}$. Its restriction to
 $ D^{p,0}$ co\"{\i}ncides with  $-d$, where $d$ is the operator of exterior differentiation
of   $p$-forms. There exists a second natural extension of $d$ namely the adjoint
second Bianchi map $\tilde D$. It sends $ D^{p,q}$ into $ D^{p,q+1}$.\\ 
The operators  $\delta=c\tilde D+\tilde D c$  and $\tilde \delta=cD+Dc$ generalize the classical 
$\delta$ operator on differential forms. Furthermore, they are respectively the formal adjoints (with respect to the integral scalar product on a compact manifold) of 
the operators $D$ and $\tilde D$.\\
The operator $D\tilde D+\tilde D D$ sends a  $(p,q)$-double form to
a  $(p+1,q+1)$-double form and its restriction to functions ((0,0)-double forms) is twice the usual
Hessian:
\begin{equation}
[D\tilde D+\tilde D D](f)=2\hhh(f).
\end{equation}
Similarly, for $p,q\geq 1$,  the operator $\delta\tilde\delta+\tilde\delta\delta$ sends
a  $(p,q)$-double form to
a  $(p-1,q-1)$-double form and satisfies
\begin{equation}\label{G}
\delta\tilde\delta+\tilde\delta\delta=(-1)^{(p+q)(n-p-q)}*(D\tilde D+\tilde D D)*.
\end{equation}
Furthermore, with respect to the integral scalar product on a compact manifold, the operator
$\delta\tilde\delta+\tilde\delta\delta$ is the formal adjoint of the operator  
$(-1)^{(p+q)}(D\tilde D+\tilde D D)$.\\

\subsection{The first variation formula}
Let  $F$ be a local variation of $M$, that is a smooth map
$$F:M\times (-\epsilon, \epsilon ) \rightarrow \tilde M,$$
such that $F(x,0)=x$ for all $x\in M$ and with  compact support ${\rm supp}F$, where
\begin{equation*}
{\rm supp}F=\overline{\{x\in M\,\, \exists t \in (-\epsilon, \epsilon )\, \, F(x,t)\not= x \}}.
\end{equation*} 
The implicit function theorem implies that there exists $\epsilon >0$ such that for all
$t$ with $|t|<\epsilon$, the map
$\phi_t=F(.,t):M\rightarrow \tilde M$ is a diffeomorphism onto a submanifold $M_t$ of $\tilde M$.
\\
Denote by $g_t$ the pull back via $\phi_t$ of the induced metric on $M_t$ from
 $(\tilde M,\tilde g)$, precisely $g_t=\phi_t^*(\tilde g)$. Note that $g_1=g$.\\
\begin{lemma}\label{lemma1} If $\xi={\frac{d}{dt}}_{|t=0}\phi_t$ denotes the variation vector 
field relative to $F$, then the first variation of $g_t$ is given by
\begin{equation}\label{F}
h(u,v)={\frac{d}{dt}}_{|t=0}g_t(u,v)=2B_N(u,v)+A_{\xi^T}(u,v).
\end{equation}
Where $N,\xi^T$ are respectively the normal and tangent components
 of the vector field
$\xi$,  $A_{\xi^T}(u,v)=g(\nabla_u\xi^{T},v)+g(\nabla_v\xi^{T},u)$ is like $B_N$ a
 symmetric bilinear form
and $\nabla$ denotes the Levi-Civita connection of $(M,g)$.
\end{lemma}
{\sl Proof}. Using for example  coordinate vector fields, one can prove without
 difficulties that
\begin{equation*}
h(u,v)=\tilde g(\tilde \nabla_u\xi,v)+\tilde g(\tilde \nabla_v\xi,u).
\end{equation*}
Where $\tilde \nabla$ denotes the Levi-Civita connection of $\tilde M$. So if $\xi= N+\xi^T$,
where $N=\xi^\bot$, then $\tilde g(\tilde \nabla_u\xi^T,v)=g(\nabla_u\xi^{T},v)$ and
$\tilde g(\tilde \nabla_uN,v)=B_N(u,v)$.\ppp \\

\begin{lemma} Let $X$ be a tangent vector field to $M$, $T$  a  symmetric
$(1,1)$-double form on $M$ and let $\alpha=T(X,.)$ be the corresponding 1-form, then
\begin{equation}\label{F+}
\langle T,A_X\rangle =\delta \alpha -\delta T.
\end{equation}
Where $A_X(u,v)=g(\nabla_uX,v)+g(\nabla_vX,u)$.
\end{lemma}
{\sl Proof.} Let $\{e_i\}$ be local orthonormal vector fields around $m\in M$ which 
diagonalize $T$ at $m$ and such that $\nabla_{e_i}=0$, then at $m$ we have
\begin{equation*}
\langle T,A_X\rangle=2\sum_iT(\nabla_{e_i}X,e_i)=\sum_i \left\{ \nabla_{e_i}(\alpha(e_i))-(\nabla_{e_i}T)(X,e_i)
\right\}.
\end{equation*}
\ppp \\
We are now ready to state and prove the first variation formula:
\begin{theorem}\label{Main}
Let $M$ be a  submanifold  of the Riemannian manifold $(\tilde M,\tilde g)$. Let 
$\xi={\frac{d}{dt}}_{|t=0}\phi_t$ denotes the variation vector 
field relative to a local variation $F$  of $M$ with compact support as above.\\
\begin{enumerate}
\item
 If 
$H_{2k}(t)=\int_M h_{2k}(g_t)\mu_{g_t}$ denotes the total $(2k)$-th Gauss-Bonnet curvature of $\phi_t(M)$, 
where
 $\mu_{g_t}$ denotes the corresponding Riemannian volume element, then 
\begin{equation}
H_{2k}'(0)=\int_M h_{2k+1}(\xi^\bot)\mu_g.
\end{equation}
\item
The submanifold $M$ is a critical point for the total 
$(2k)$-th Gauss-Bonnet
curvature function for all local variations of $M$ if and only if the $(2k+1)$-Gauss-Bonnet
curvature $h_{2k+1}(N)$ of $M$  vanishes for all normal directions $N$.
\end{enumerate}

\end{theorem}
{\sl Proof}. 
Since by construction  $h$ vanishes outside the compact subset ${\rm supp} F$, 
lemma \ref{lemmatpq1} below
 and Gauss' theorem
imply that
\begin{equation*}
\begin{split}
H_{2k}'.h&=\int_M\biggl(h_{2k}'.h+\frac{h_{2k}}{2}
{\rm tr}_gh\biggr)\mu_g\\
&=-\frac{1}{2}<\frac{c^{2k-1}}{ (2k-1)!}R^k,h>+\frac{h_{2k}}{2}<g,h>\\
&=\frac{1}{2}<h_{2k}g-\frac{c^{2k-1}}{(2k-1)!}R^k,h>\\
&=\frac{1}{2}\langle T_{2k},h\rangle.
\end{split}
\end{equation*}
Where $\langle ,\rangle$ is the integral scalar product. Next using (\ref{F}), (\ref{F+}), 
Gauss's theorem and the fact that  Einstein Lovelock tensors are divergence free
\cite{Labbivariation} we get
\begin{equation*}
\frac{1}{2}\langle T_{2k},h\rangle=\langle T_{2k},B_N\rangle .
\end{equation*}
Consequently, equation (\ref{E}) shows that $H_{2k}'.h=h_{2k+1}(N)$ as desired. \ppp

\begin{lemma}[\cite{Labbidoubleforms}]\label{lemmatpq1}
The directional derivative of $h_{2k}$ at $g$ in a given
direction  $h$  is given by
\begin{equation*}
h_{2k}'h=\frac{-1}{2}\langle \frac{c^{2k-1}}{(2k-1)!}R^k,h\rangle +{\rm div}W_1+{\rm div}W_2.
\end{equation*}  
Where $W_1$ (resp. $W_2$) is the tangent vector field over $M$ corresponding to the $1$-form
$\delta \biggl( *( \frac{kg^{n-2k}}{4(n-2k)!}R^{k-1}h)\biggr)$
(resp. $\tilde\delta \biggl( *( \frac{kg^{n-2k}}{4(n-2k)!}R^{k-1}h)\biggr)$).
\end{lemma}

\subsection{$(2k)$-Minimal submanifolds}
With respect to the previous variational formula and by analogy to the case of 
usual minimal submanifolds we set the following definition:
\begin{definition} For $0\leq 2k\leq n$,
An $n$-submanifold $M$ of the Riemannian manifold $(\tilde M,\tilde g)$ is said to be $(2k)$-minimal if
\begin{equation*}
h_{2k+1}\equiv 0.
\end{equation*}
\end{definition}
\par\medskip\noindent
In the extreme cases: the $0$-minimal submanifolds are nothing but the usual minimal 
submanifolds. And if $n$ is even, every submanifold
is $n$-minimal (the condition is empty). \\
We provide below examples of intermediate minimal submanifolds:
\begin{enumerate}
\item A flat submanifold is always $(2k)$-minimal for all $k> 0$.  In fact 
$R\equiv 0\Rightarrow
h_{2k+1}\equiv 0$. This shows that $(2k)$-minimal does not imply the usual minimality
condition.\\
\item A totally geodesic submanifold is always $(2k)$-minimal for all $k\geq 0$.  In fact $B\equiv 0\Rightarrow
h_{2k+1}\equiv 0$.\\
\item A submanifold with constant curvature $\lambda \not= 0$ is $(2k)$-minimal if and only if it is
minimal in the usual sens. In fact, in this case $R=\frac{\lambda}{2} g^2$, and therefore 
$$h_{2k+1}(N)=*(\frac{g^{n-2k-1}}{(n-2k-1)!}2^{-k}\lambda^kg^{2k}B_N)=
\frac{(n-1)!\lambda^k}{(n-2k-1)!2^k}cB_N.$$
\\
\item If $M$ is a hypersurface of the Euclidean space then $(2k)$-minimality coincides with
Reilly's $(2k)$-minimality \cite{Reilly}. In fact, as we have seen above, in this case $h_{2k+1}$
coincides, up to a constant, with the $(2k+1)$-mean curvature of the hypersurface.\\
In particular, there are no $(2k)$-minimal compact hypersurfaces of the Euclidean
space. For, it is standard that there always exists one point on the hypersurface where
all the mean curvatures are positive and in particular not equal to zero.\\
\item If $M$ is a hypersurface of a space form $(\tilde M,\tilde g)$ of constant $\lambda$ then
$M$ is $(2k)$-minimal if and only if 
$$\sum_{i=0}^k\frac{k!(2k-2i+1)!(n-2k-1+2i)!\lambda^i}{i!(k-i)!(n-2k-1)!2^k}s_{2k-2i+1}=0.$$
Where $s_j$ denotes the symmetric functions in the eigenvalues of the shape
 operator of the hypersurface. This fact can be proved easily after using the Gauss equation
$R=\lambda/2 g^2+1/2B^2$ and the binomial theorem. Notice the difference with Reilly's
 $r$-minimality \cite{Reilly}.
\item Any complex submanifold $M$ of a Kahlerian manifold $(\tilde M,\tilde g)$ is
 $(2k)$-minimal for any $k$. \\
To prove this fact,  note first that $M$ is then itself a Kahlerian manifold.
The Riemann curvature tensor of $M$ is therefore invariant 
under the complex structure $J$, that is $R(J.,J.,J.,J.)=R(.,.,.,.)$.\\
 It can be shown without difficulties that the 
 same property is therefore true for the tensors $R^k$ and consequently for the contraction
 $c^{2k-1}R^k$. Consequently, 
the Einstein-Lovelock tensors 
$T_{2k}=h_{2k}g-\frac{c^{2k-1}R^k}{(2k-1)!}$ are $J$-invariant, that is  $T_{2k}(J.,J.)=T_{2k}(.,.)$.\\
 On the other hand, the 
second fundamental form of $M$
satisfies $B(Jx,y)=B(x,Jy)$. It is then straightforward that $h_{2k+1}(N)=
\langle T_{2k},B_N\rangle \equiv 0$.\\
\item Suppose the submanifold $(M,g)$ is $(2k)$-Einstein,
that is $T_{2k}=\lambda g$, then:
\begin{itemize}
\item If $\lambda=0$ then $M$ is $(2k)$-minimal.
\item If $\lambda\not=0$ then $M$ is $(2k)$-minimal if and only if it is minimal 
in the usual sens.
\end{itemize}
In particular,  since irreducible isotropy homogeneous spaces are $(2k)$-Einstein for all $k>0$ 
\cite{Labbivariation}, then they admit a $(2k)$-minimal immersion in a sphere for all $k>0$.
\end{enumerate}
{\sc Remark.} Since the considerations of theorem \ref{Main} are  of local nature,  the theorem
 remains then  true for immersed
submanifolds $M\subset \tilde M$. Therefore it makes sens to consider also
 $(2k)$-minimal immersed submanifolds.\\

\subsection{Generalized Laplace Operators}
Let $(M,g)$ be a Riemannian manifold and $f$ be a smooth function defined on $M$.
Recall that the usual Laplacian of $f$ is given by
$$\Delta f=-c{\rm Hess}\, ( f)=-\langle g, \hhh (f)\rangle  $$
where  $\hhh (f)$ denotes the Hessian of the function $f$.\\
Instead of just taking  the trace of the Hessian, one can takes the determinant
 (Monge-Amp\`ere operator) and more generally the elementary symmetric functions
 in the eigenvalues of $\hhh (f)$. These generalized operators
 appear naturally
 in the context of the $\sigma_k$-Yamabe problem, 
where they are denoted $\sigma_k(\hhh(f))$, see \cite{Viacklovsky} 
and the references therein.
 Precisely, they are up to a constant equal to
 $$c^{2k}\hhh^k(f)=\langle \hhh^k(f),g^k\rangle. $$
 In our context of generalized minimal submanifolds, another natural generalization
 of the Laplace operator appears naturally. Precisely, we set
 \begin{definition}
 Let $f$ be a smooth function on $(M,g)$. We define the $\ell_{2k}$-Laplacian
 operator
 of $(M,g)$ as
 \begin{equation}
 \ell_{2k}(f)=-\langle T_{2k},\hhh(f)\rangle.
 \end{equation}
 Where $T_{2k}$ denotes the $(2k)$-th Einstein-Lovelock tensor of $(M,g)$ and
 $0\leq 2k<n$.\\
We shall say that the function $f$ is $\ell_{2k}$-harmonic if $\ell_{2k}(f)=0.$
 \end{definition}
For $k=0$ we have $T_0=g$ and then $\ell_0=\Delta$ is the usual Laplacian. Furthermore,
if $(M,g)$ is $(2k)$-Einstein, that is $T_{2k}=\lambda g$, then
$\ell(f)=\lambda \Delta(f)$ is, up to a constant, the usual Laplacian.\\
In particular, the later property holds for manifolds with constant curvature and
for isotropy irreducible homogeneous manifolds \cite{Labbivariation}.\\
We prove below some properties of these operators.
\begin{proposition}
Let $D,\tilde D, \delta$ and $\tilde \delta$ be as in section \ref{Diff}. 
The $\ell_{2k}$-Laplacian can be written in any one of the following
equivalent forms: 
\begin{equation*}
\begin{split}
\ell_{2k}(f)&=*\left\{ \frac{g^{n-2k-1}}{(n-2k-1)!}R^k\hhh(f)\right\}\\
&= *\left\{ [D\tilde D+\tilde DD](\frac{g^{n-2k-1}}{2(n-2k-1)!}fR^k)\right\}\\
&=\frac{1}{2}[\delta\tilde\delta+\tilde\delta\delta](fT_{2k}).\\
\end{split}
\end{equation*}
In particular, $\ell_{2k}(f)$ is a divergence and therefore
$\int_M\ell_{2k}(f)dv\equiv 0$ if $M$ is compact.
\end{proposition}
{\sl Proof.} The first formula is a direct consequence of the definition of 
$T_{2k}$ and formula (\ref{B}). The second one results from the fact that the metric $g$
and the Riemann tensor are in  $\ker (D \cap \tilde D)$, and the fact that the former
is closed under the exterior product of double forms  \cite{Kulk}.
The last formula results from (\ref{B+}) and (\ref{G}).\\
\begin{proposition}
If for some $k$ with $0\leq 2k<n$,
 the Einstein-Lovelock tensor $T_{2k}$ is positive definite (or negative definite), then
the operator $\ell_{2k}$ is elliptic.
\end{proposition}
{\sl Proof}. Recall that in local coordinates we have
$$ \hhh (f)=(\frac{\partial^2f}{\partial x^i\partial x^j}-\frac{\partial f}{\partial x^k}\Gamma_{ij}^k)dx^i\otimes dx^j.$$
Therefore,
$$ \ell_{2k}(f)=\sum_{i,j}T_{2k}(\frac{\partial }{\partial x^i},\frac{\partial }{\partial x^j})
(\frac{\partial^2f}{\partial x^i\partial x^j}-\frac{\partial f}{\partial x^k}\Gamma_{ij}^k).$$
\ppp
\\
\begin{proposition}
If $M$ is compact, then the operator $\ell_{2k}$ is self adjoint with respect to the integral scalar product, that is for
arbitrary smooth functions $u,v$ on $M$ we have
\begin{equation*}
\langle u,\ell_{2k}(v)\rangle=\langle v,\ell_{2k}(u)\rangle.
\end{equation*}
Furthermore, If $T_{2k}$ is positive definite (resp. negative definite) then $\ell_{2k}$ is
positive definite (resp. negative definite). More precisely we have the following integral formula
\begin{equation*}
 \langle \ell_{2k}(f), f\rangle=T_{2k}(df^{\sharp},df^{\sharp}).
\end{equation*}
\end{proposition}
{\sl Proof.} First note that
\begin{equation*}
\begin{split}
\hhh(uv)&=\frac{1}{2}[D\tilde D+\tilde DD](uv)\\
&=v\hhh(u)+u\hhh(v)+(Du\tilde D v+\tilde D u Dv)\\
&=v\hhh(u)-u\hhh(v)+D(u\tilde D v)+\tilde D(uDv).
\end{split}
\end{equation*}
Since, with respect to the integral scalar product, the formal adjoints of $D$ and $\tilde D$
are divergences, namely $\delta$ and $\tilde \delta$, see section \ref{Diff},
 and since the Einstein Lovelock tensors are divergence free \cite{Labbivariation}, it results from the previous formula that 
\begin{equation*}
\begin{split}
0=-\int_M \ell_{2k}(uv){\rm dvol}&=\langle T_{2k}, \hhh(uv)\rangle\\
&=\langle T_{2k}, v\hhh(u)\rangle-\langle T_{2k}, u\hhh(v)\rangle+0\\
&=-\langle v,\ell_{2k}(u)\rangle+\langle u,\ell_{2k}(v)\rangle .\\
\end{split}
\end{equation*}
Next, it results from the formula in the third line of this  proof, after taking $u=v=f$, that 
\begin{equation*}
0=-2\langle \ell_{2k}(f), f\rangle +\langle Df\tilde D f+\tilde D f Df ,T_{2k} \rangle .\\
\end{equation*}
On the other hand, at each point of $M$ we have
\begin{equation*}
\begin{split}
\langle Df\tilde D f+\tilde D f Df ,T_{2k} \rangle &=2\sum_{i,j}T_{2k}(e_i,e_j)df(e_i)df(e_j)\\
=& 2\sum_{i,j}T_{2k}(df(e_i)e_i,df(e_j)e_j)=2T_{2k}(df^{\sharp},df^{\sharp}).
\end{split}
\end{equation*}
This completes the proof. \ppp \\
The proof of the following corollary is straightforward:
\begin{corollary}\label{harmonic}
Let $(M,g)$ be a compact manifold of positive definite (or negative definite) Einstein-Lovelock
tensor $T_{2k}$ then every  smooth and $\ell_{2k}$-harmonic function on $M$ is constant.
\end{corollary}
\subsection{$(2k)$-minimal submanifolds in Euclidean space and in the spheres}
We suppose now that $\tilde M=\Re^{n+p}$ is the Euclidean space. For  $v\in \Re^{n+p}$, we
define the coordinate  function $f_v:M\rightarrow \Re$ by $f_v(m)=\langle v, m\rangle $. It is 
not difficult to see that the Hessian of $f_v$ is nothing but the second fundamental form
of $M$ in the direction of the normal component of $v$, that is
$$\hhh(f_v)(x,y)=-\langle B(x,y),v\rangle .$$
In particular if $v$ is normal to $M$, we get
\begin{equation}
\ell_{2k}(f_v)=-\langle T_{2k},\hhh(f_v)\rangle =\langle T_{2k},B_v\rangle=h_{2k+1}(v).
\end{equation}
We have therefore proved the following result:
\begin{proposition}
A submanifold $M$ of the Euclidean space is $(2k)$-minimal if and only if the coordinate
functions restricted to $M$ are $\ell_{2k}$-harmonic functions on $M$.
\end{proposition}
\par\medskip\noindent
Let now $F:M^n\rightarrow \Re^{n+p}$ be an isometric immersion  and let 
$F_i(x)=\langle F(x),e_i\rangle $
be the $i$-th component of $F$ in $\Re^{n+p}$. If we look to $F$ as a function
from $F(M)$ to $\Re^{n+p}$and since the results are local, me may assume $F(M)$ a submanifold and then
$\hhh(F_i)=B_i$. In particular, 
if $h_{2k+1}$ is the Gauss-Bonnet curvature of 
$F(M)$ in $\Re^{n+p}$ then we  we have (componentwise)
\begin{equation}
\ell_{2k}(F)=h_{2k+1}.
\end{equation}
In particular,
$F$ is $(2k)$-minimal if and only if $\ell_{2k}(F_i)=0$ for all $i$.\\
Using corollary \ref{harmonic} and the previous remark we immediately obtain
\begin{corollary}
Let $0\leq 2k<n$ and let $(M,g)$ be a compact Riemannian  $n$-manifold with 
 positive definite (or negative definite) Einstein-Lovelock tensor $T_{2k}$. Then there is
no non trivial isometric $(2k)$-minimal immersion of $M$ in the Euclidean space. \\
In particular, there are no non trivial compact $2$-minimal submanifolds in the Euclidean space
with positive definite (or negative definite) Einstein tensor. In other words, if
$M$ is a compact submanifold of the Euclidean space with positive definite (or negative definite)
Einstein tensor, then there exist variations of $M$ which increase the total scalar curvature
and others which decrease it.\\
\end{corollary}
Note that the condition of positive (or negative) definiteness of $T_{2k}$ in the previous corollary is
necessary, as the flat torus admits non trivial $(2k)$-minimal isometric immersions in the Euclidean space.
\par\medskip\noindent
 Finally, we  prove the following about $(2k)$-minimal immersions in spheres.
\begin{proposition}
Let $F: M^n\rightarrow S^{n+p}\subset \Re^{n+p+1}$ be an isometric immersion. Then $F$
is $(2k)$-minimal into $S^{n+p}$ if and only if there is a smooth function $\phi:M\rightarrow \Re$
such that $\ell_{2k}F=\phi F$ (componentwise).
\end{proposition}
{\sl Proof.} It is not difficult to see that the two second fundamental forms of the 
submanifold in the
sphere and in the Euclidean space coincide in  normal directions that are tangent
to the sphere. Consequently, it results from its definition that the $(2k+1)$th Gauss-Bonnet
 curvature $h_{2k+1}$ of $F(M)$ in $S^{n+p}$
coincides with the restriction to $TS^{n+p}$ of the $(2k+1)$th 
Gauss-Bonnet curvature of $F(M)$ in $\Re^{n+p+1}$. The later being equal to $\ell_{2k}(F)$
as above, then $h_{2k+1}\equiv 0$ if and only if the components of $\ell_{2k}(F))$ vanishes for the
directions tangent to the  sphere. 
But since $F$ takes its values in the sphere then 
$F(x)$ are normal vectors to $S^{n+p}$ therefore $F$ is a $(2k)$-minimal if and only if
$\ell_{2k}(F)=\phi F.$\ppp

\vspace{2cm}
\noindent
Labbi M.-L.\\
  Department of Mathematics,\\
 College of Science, University of Bahrain,\\
  P. O. Box 32038 Bahrain.\\
  E-mail: labbi@sci.uob.bh

\end{document}